# Dynamic Network Models

Benjamin Armbruster       John Gunnar Carlsson

March 12, 2011

**Abstract**

We analyze random networks that change over time. First we analyze a dynamic Erdős-Renyi model, whose edges change over time. We describe its stationary distribution, its convergence thereto, and the SI contact process on the network, which has relevance for connectivity and the spread of infections. Second, we analyze the effect of node turnover, when nodes enter and leave the network, which has relevance for network models incorporating births, deaths, aging, and other demographic factors.

## 1  Introduction

We present here a theoretical analysis of some dynamic models of undirected random networks. We have two major contributions in this paper. First, we analyze a dynamic Erdős-Renyi model, in particular its stationary distribution, its convergence thereto, and the SI contact process on the network, which has relevance for connectivity and the spread of infections. Second, we analyze the effect of node turnover (when nodes enter and leave the network) on dynamic Erdős-Renyi networks and preferential attachment networks, which has relevance for network models incorporating births, deaths, aging, and other demographic factors. By analyzing a dynamic network model, this work brings together two streams of literature, one in the social sciences that focuses on the statistical problem of fitting dynamic network models to data and one in mathematics and physics that focuses on analyzing the properties of static networks and processes on them. In the next section we present an overview of the existing literature. In section 3 we analyze the dynamic Erdős-Renyi $G(n,p)$ model, and in section 4 we present the model of node turn-over.

## 2 Literature

### 2.1 Social Sciences

Work on dynamic random graphs in the social sciences dates back to 1977 [16, 31]. In [16], Holland and Leinhardt describe a general framework for modeling the evolution of a graph as a continuous time Markov chain where one edge changes at a time. They focus on directed networks and a simplification of their simplest model is essentially the dynamic Erdős-Renyi model we discuss in the next section, which is sometimes referred to as the "binomial" model. The driving problem in this stream of literature is to develop and parameterize models that can describe and explain longitudinal social network data (i.e., data at multiple points in time of the network, the characteristics of both the nodes and links). Thus this literature has focused its attention on developing models with sophisticated dependence structure and the statistical tools necessary to fit them to data. See [28] for a current snapshot of the field and [27] for a more thorough overview. Our work complements this stream of research by using similar models but instead of examining statistical questions we focus on analytical aspects such as how long will it take an infection to spread through the entire network and how this depends on the network properties.

**Node Turnover** In section 4 focuses specifically on node turnover. In the social science literature, [17] has looked at the simulation and statistical aspects of a similar problem. They extend the popular and flexible stochastic actor-oriented framework to incorporate changing node sets. Our contribution is to analyze the connection between the turnover process and the structural properties of the network, in particular the degree distribution.

### 2.2 Physics and Mathematics

**Analysis of static random graphs** There is a large literature in physics and math analyzing static random graphs and processes on them. Two of the most popular focuses of this line of work are Erdős-Renyi networks and preferential attachment models. Rigorous results for Erdős-Renyi networks date back to the 1959 paper by Gilbert [13] while preferential attachment models were popularized by Barabási and Alberts in 1999 [1]. In addition rigorous probabilistic analyses, approximations are popular in this stream of literature as they greatly increase the range of models that can be analyzed. Of particular note are mean-field approximations that describe the network as a set of differential equations tracking over time the number of nodes with certain characteristics such as having a certain degree, being infected, or a combination of the two. The books [9] and [3] provide an overview of this stream of literature with the former more mathematically rigorous than the latter.

**SI contact process on Erdős-Renyi networks** Much of the popularity of Erdős-Renyi networks is due to ability to prove many results about the size of their connected components [10]. These results can be interpreted as the final epidemic size of an SI contact process. The mean-field approximation of the SI contact process on Erdős-Renyi networks is simply $\dot{x} = \beta d x(1-x)$ where $x$ is the fraction of infected nodes, $\beta$ is the rate of infection across each edge, and $d$ is the average degree. However, we are not aware other, more rigorous work characterizing the SI contact process on Erdős-Renyi networks. The related SIS contact process has been studied using a mean-field approximation by Pastor-Satorras and Vespignani [21] on both preferential attachment models and small-world networks, where the analysis for small-world networks applies just as well to Erdős-Renyi networks.

## 2.3 Analysis of Dynamic Networks

**Dynamic Erdős-Renyi model** The only previous analysis of the dynamic Erdős-Renyi model we are aware of is [12]. In 2009, [12] laments: " 'What is the dynamic analog of the Erdős-Rnyi random graph model $G(n, p)$?' No one appears to have defined such a thing[.]" The authors then introduce a dynamic Erdős-Renyi model in the context of a counter-terrorism application. The model provides the authors a more realistic, time-correlated network noise structure for tracking threatening behavior in networks. Their analysis limits itself to noting that the graph at any fixed time is an Erdős-Rnyi graph and likewise for the union or intersection of the graph over an interval of time. We significantly expand upon this analysis to the extent that we suggest that the analytical tractability of this model makes it the canonical dynamic network model.

**Dynamic Random Geometric Graphs** A similar analysis to ours for dynamic Erdős-Renyi graphs has been made for dynamic random *geometric* graphs [22], which have important applica-tions for mobile ad-hoc networks. These are graphs where the vertices (the communication nodes) are embedded in $\mathbb{R}^d$ and their locations follow random walks. Edges exist between any two vertices whose distance is below a threshold (i.e., the nodes can communicate). The authors determine asymptotics for the length of the periods where the graph is connected and the length of time between those periods.

**Preferential Attachment** One can also view the preferential attachment model and its variations as dynamic random graphs because in the preferential attachment model nodes and edges are added over time. See [4, 5] for two recent surveys. These models differ in three significant ways from the dynamic Erdős-Renyi model we study: they are in discrete time; the number of nodes and edges grows without bound (even in those models that allow for node and edge removal); and most importantly these models focus on the (asymptotic) structure of the resulting (infinite) graph and view the dynamics as merely a means of generating that structure. We see the dynamic Erdős-Renyi model, with its finite size, as being uniquely suited for focusing on the *dynamics* of the model rather than the structure of the network.

**Pair Model** The above exhausts the mathematically rigorous analysis of dynamic random graphs we were able to find. We now survey the work using some approximations to analyze dynamic random graphs. In addition to the flexible mean-field approximations, a pair model has been to analyze dynamic networks. In fact, the first analytical approximation for disease spread on a dynamic network is the pair-model [8]. This is an ODE model that goes beyond the standard mean-field models by having dependent variables for the number of single infected males, females, and pairs (one each for the kind where only the male is infected, only the female, and both partners). It approximates a dynamic network of constant size by allowing partnerships to form and dissolve. More recently, this approach has been used in the physics literature to approximate stylized disease spread in a variety dynamic networks [15, 23, 30].

**Rewiring** While the network dynamics studied with pair-models and mean-field approximations include models with edges appearing and disappearing (the previously mentioned [8, 23]) or with two pairs of nodes swapping edges (the previously mentioned [30]), they most often use a "rewiring" network dynamic. Rewiring involves changing one endpoint of an edge at a time and is motivated by the process generating the long-range links of a small-world network. While disease spread on a dynamic network with rewiring has sometimes been analyzed using a pair model approximation [15], it is most often studied by a combination of simulation and the mean-field approximation. The models range from ones where the detachment and attachment probabilities of rewiring depend on node degree [20, 25] and centrality [11] to models where the rewiring is based on the state of the node (e.g., preferentially breaking links between susceptible and infected nodes) [14, 15, 24, 26, 29]. These latter models are sometimes called "adaptive" or co-evolutionary" networks. Surprisingly, we do not know of any work analyzing the dynamic Erdős-Renyi model using a mean-field approximation or a pair model despite it being amenable to these approaches and it being simpler than most of the models mentioned previously (with possibly the exception of [8]). Instead of using an ODE approximation that elides the network structure (such as a mean-field approximation or a pair model), we instead rigorously analyzes the dynamic Erdős-Renyi model.

**Node-turnover** Our model of node turn-over finds a power-law degree distribution and is an extension of either the dynamic Erdős-Renyi model or the preferential attachment model. Our model has similarities to variations of the preferential attachment model which also allow for node removals. The most similar are the *Cooper-Frieze-Vera* [7] and *Chung-Lu* [6] models. As in our model (the one extending the preferential attachment model), these models add edges to new nodes based on a preferential attachment rule and delete nodes uniformly at random. (These models also allow for some events other than node addition or removal.) Nevertheless, these models, like all the variations of preferential attachment models we are aware of and all the dynamic models in the survey [5], differ fundamentally from our model because their number of nodes goes to infinity, while in ours the number remains roughly constant, which we believe is helpful for many modeling applications. The only analysis of network models with node turnover whose number of nodes does not go to infinity is recent work in 2010 by Christel Kamp [18, 19] who simulates a mean-field approximation of a fatal disease spreading on a network. We do not look at disease spread with our model of node turn-over and focus more on deriving analytical results.

# 3 Dynamic $G(n,p)$ networks

**Definition and Stationary Distribution** In our model we suppose that we have $n$ nodes and we let $e_{ij}(t)$ equal 1 if there is an edge between nodes $i$ and $j$ and 0 otherwise.

We first present a model extending the simplest random networks, the Erdős-Renyi $G(n,p)$ graphs. As in those graphs, we assume the state of each potential edge, $e_{ij}$, is independent, but instead of each edge being a Bernoulli random variable, $e_{ij}(t)$ is a telegraph process, a 2-state CTMC with rate $\lambda$ from $0 \to 1$ and rate $\mu$ from $1 \to 0$. The rate at which this CTMC completes a cycle is $\alpha := \mu\lambda/(\mu+\lambda)$ and its stationary distribution is $\bar{e}_{ij} \sim$ Bernoulli$(p)$ where $p := \lambda/(\mu+\lambda)$. Thus the stationary distribution of the graph is a $G(n,p)$ graph. We will sometimes denote this dynamic graph as $G(n,\lambda,\mu)$ or $G(n,p,\alpha)$. Obviously, $\alpha = 0$ is the case of a static $G(n,p)$ graph. Note that $\mu = \alpha/p$ and $\lambda = \alpha/(1-p)$. We will let $d := (n-1)p$ denote the expected degree of any node in steady-state.

**Convergence to Stationary Distribution** The distribution of the edge state approaches the stationary distribution exponentially at a rate $\lambda + \mu = \frac{\alpha}{p(1-p)}$. Specifically, $\Pr[e_{ij}(t) = 1] - p = (\Pr[e_{ij}(0) = 1] - p)e^{-(\lambda+\mu)t}$. Now let's consider the distribution of $k$ edges, $E_k(t) := (e_1(t), \ldots, e_k(t))$, converging to their stationary distribution, $\bar{E}_k := (\bar{e}_1, \ldots, \bar{e}_k)$. We will consider the case with the slowest convergence where $p \leq 1/2$ and for every edge $e_i(0) = 1$. This will simplify the calculations without affecting the rates of convergence. A sensible measure of convergence is the total-variation distance and the related mixing time. The total-variation between any two measures $\xi$ and $\nu$ is the largest difference across all events: $TV(\xi, \nu) := \sup_E |\xi(E) - \nu(E)|$. The mixing time $\tau$ is the time until $TV(E_k(t), \bar{E}_k) = 1/4$. The mixing times below are fast. The proof is in the appendix.

**Theorem 1.** *Suppose that $k \to \infty$. If $p$ constant, then $\tau \sim \frac{1}{2(\lambda+\mu)} \log k$ and if $p := c/k$ for some $c$, then $\tau \sim \frac{c \log k}{\alpha k}$.*

**SI Contact Process** Now let's look at the SI contact process with a rate of infection $\beta$ across each edge. This is sometimes called the "contagion" or "diffusion" process. Two important special case of this process are

1. $\beta = \infty$. This case is a dynamic analog of graph connectivity. We will discuss this in detail subsequently.

2. $\alpha = \infty$. This case is equivalent to the (static) case with $p' = 1$ and $\beta' = \beta p$.

For the SI process we want to look at the counting process $X(t)$ counting the number of infected nodes. We're interested in both the distribution of $X(t)$ and it's hitting times, $\tau_k := \inf\{t : X(t) \geq k\}$. There is a simple time-scaling: $X(rt)$ and $\tau_k/r$ are equivalent to $X'(t)$ and $\tau'_k$ respectively, with $\alpha' = r\alpha$ and $\beta' = r\beta$. We assume that at time 0 we're in steady-state, we start with a $G(n,p)$ graph.

## 3.1 Case $\beta = \infty$

This case is a dynamic analog of graph connectivity. Recall the following theorem about connectivity of Erdős-Renyi graphs.

**Theorem 2** (Theorem 2.8.1 in [9]). *As $n \to \infty$, the probability $G(n,p)$ with $p = a \log n / n$ is connected goes to 0 if $a < 1$ and to 1 if $a > 1$.*

Thus if $p > 0$, then $\Pr[\tau_n = 0] \to 1$ as $n \to \infty$. However for $p = 0$ we can obtain an asymptotic lower bound on $\tau_n$. This graph has $\lambda = \alpha$ and $\mu = \infty$. By ignoring the removal of edges we obtain a lower bound: setting $\mu = 0$ instead of $\mu = \infty$. Thus the following theorem implies that $\tau \succeq \log n/(\alpha n)$ in probability as $n \to \infty$ when $p = 0$.

**Theorem 3.** *For a $G(n, \lambda, \mu)$ graph with $\mu = 0$, $\tau_n \sim \log n/(\lambda n)$ in probability as $n \to \infty$.*

**Upper bound for $\tau_n$** We're going to use coupling and create a Markov birth process $X'(t)$ such that $X'(t) \leq X(t)$. We assume $X'(0) = 0$, that initially there are no edges. Then we assume only one edge is added at a time, $X'$ transitions from $k$ to $k+1$ at rate $\lambda k(n-k)$. Thus,

$$\mathbb{E}[\tau_k] \leq \sum_{i=1}^{k-1} \frac{1}{\lambda i(n-i)} = \frac{1}{\lambda n}(H_{k-1} + H_{n-1} - H_{n-k}) \tag{1}$$

where $H_k = \sum_{i=1}^{k} 1/i$ is the $k$th harmonic number. Since $H_k \leq \log k + 1$,

$$\mathbb{E}[\tau_k] \leq \frac{1}{\lambda n}\left(2 + \log \frac{nk}{n-k}\right) \tag{2}$$

and $\mathbb{E}[\tau_n] \leq 2(1 + \log(n-1))/(\lambda n)$. Asymptotically, $H_n \sim \log n$ and thus, $\mathbb{E}[\tau_n] \preceq 2\log(n)/(\lambda n)$. Note that this differs by a factor of two from the above lower bound for the case of $p = 0$.

## 3.2 Case $\beta < \infty$

For the case of $\beta < \infty$ we resort to asymptotics as $n \to \infty$. Proofs are in the appendix.

We start by deriving an upper bound on the time to infect a node, $\mathbb{E}[\tau_{m+1} - \tau_m]$. To bound this from above we assume that all edges in the network disappear. Thus there are $N := m(n-m)$ potential edges between infected and susceptible nodes. We build a Markov chain with $N+2$ states $(0, \ldots, N$ for the number of potential edges and an absorbing state $S$ representing a new infection occurring). Let $t_k$ be the expected time to reach the absorbing state when starting in state $k$. We would like to determine $t_0$, the expected time to hitting the absorbing state $S$ when starting in state $0$, since $\mathbb{E}[\tau_{m+1} - \tau_m] \leq t_0$. The transition rates in state $k$ are as follows: since there are $k$ edges that are currently active, we're moving towards the absorbing state at a rate $k\beta$. We're moving towards the state $k-1$ at a rate $k\mu$. We're moving towards the state $k+1$ at a rate $(N-k)\lambda$.

**Lemma 4.** *For the above Markov chain, $t_0 \sim \sqrt{\pi/(2\beta\lambda N)}$ as $N \to \infty$.*

Thus
$$\mathbb{E}[\tau_k] \preceq \sum_{m=1}^{k-1} \frac{\sqrt{\pi/(2\beta\lambda)}}{\sqrt{m(n-m)}} = \sqrt{\pi/(2\beta\lambda)} \sum_{m=1}^{k-1} \frac{1}{\sqrt{m(n-m)}} \tag{3}$$

Taking the limit as $n \to \infty$, this Riemann sum becomes an integral,
$$\mathbb{E}[\tau_k] \preceq \sqrt{\pi/(2\beta\lambda)} \int_0^{k/n} \frac{1}{\sqrt{t(1-t)}} dt = \sqrt{\pi/(2\beta\lambda)} \cos^{-1}(1 - 2k/n). \tag{4}$$

Hence $\mathbb{E}[\tau_n] \preceq \sqrt{\pi^3/(2\beta\lambda)}$ and does not depend on $n$.

**Theorem 5.** $\tau_n \succeq \sqrt{2\log n/(\beta\lambda n)}$ *in probability as $n \to \infty$.*

The fact that our bounds on $\tau_n$ do not depend on $\mu$ is reflected in figure 1 where the dependence on $\mu$ is smaller than the dependence on $\lambda$ and $\beta$.

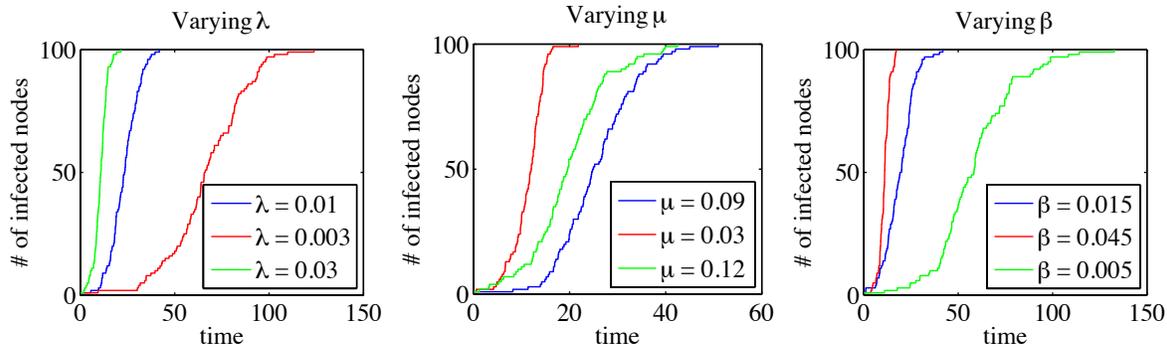

Figure 1: The number of number of infected nodes over time as we vary parameters of a dynamic Erdős-Renyi network with 100 nodes and default parameters $\mu = 0.01$, $\lambda = 0.01$, and $\beta = 0.015$.

## 4 Node Turnover

Another way a network may change is having nodes enter and leave the network. This may be due to birth, death, emigration, or immigration. When a node leaves, all the edges connected to that node also disappear. We will now let $N(t)$ be the random process counting the number of nodes in the network at time $t$.

### 4.1 Dynamic Erdős-Renyi networks

Consider a dynamic Erdős-Renyi network and suppose each node has a constant hazard rate of 1 of leaving the network. Then the network is still a Erdős-Renyi graph after the node leaves. Often it will be useful to have a network whose size is roughly constant. Thus we will assume that there is a Poisson process with rate $n$ of new nodes entering the network. When a node enters the network, it initially has no edges. Thus it is not quite an Erdős-Renyi graph after this event. However, the potential edges connected to the new nodes will soon be in steady state because their mixing time is short, $\mathcal{O}(\alpha^{-1} \log n)$.

The number of nodes in the network, $N(t)$ is a birth-death Markov process with birth rate $n$ and death rate $N(t)$. The steady-state distribution is Poisson($n$). A Poisson process with rate $n$ describes the births. Each node's life is exponentially distributed with rate 1. Thus the deaths occur in steady-state as a Poisson process with rate $n$. The entire process is time reversible.

Now consider a fixed time $t$ or a time independent of $N(\cdot)$. Due to the reversibility, each individual's age is exponentially distributed with rate 1. A node selected uniformly at random will have Poisson$(n-1)$ neighbors. The time $T$ that any *pair* of nodes exists is an exponential random variable with rate 2. Thus the probability that an edge exists is $p' := p(1 - \mathbb{E}[\exp(-(\lambda + \mu)T)]) = p(1 - (\lambda + \mu + 1)^{-1})$. Thus at time $t$ we have a $G(n,p)$ graph but with $p = p'$, a slightly smaller edge probability.

## 4.2 Preferential Attachment

Preferential attachment is often modeled as a process adding nodes one at a time to the network, leading the network to grow without bound. To model a changing network of roughly constant size, $n$, we will suppose that nodes are removed at the same rate as they are added. Our analysis won't distinguish between different ways of adding nodes at a rate of 1 per unit time (e.g., one each discrete time step; Poisson(1) each discrete time step; or with a Poisson process in continuous time). A key question is whether the degree distribution of the resulting *dynamic* network still follows a power-law. Note that the static (infinite) preferential attachment process leads to a degree distribution scaling as $k^{-3}$. We show that depending on the details of node removal we can obtain power-law degree distributions with different exponents.

Following the original preferential attachment model, we assume that each new node creates $m$ links where each link connects the new node to a random node with a probability proportional to the degree of the node. In the long-term, the total number of edges in the network is $nm$. The sum of degrees of the nodes is twice that. This will remain roughly constant. Suppose $k(t)$ is the degree of a node of age $t$. Using an argument similar to one in [2], we treat $k(t)$ as a continuous variable. Then, $k(0) = m$ and $k'(t) = mk/(2nm) = k/(2n)$. Thus, $k(t) = me^{t/(2n)}$. Suppose $F_l(\cdot)$ is the distribution of a node's lifespan. Then the density of nodes with age $t$ is $f_a(t) = (1 - F_l(t))/n$. Let $F_a(t) := \int_0^t f_a(s)ds$. Then $\Pr[k(t) > k] = (1 - F_a(2n \log(k/m))$. Differentiating,

$$\Pr[k(t) = k] = (2n/k)f_a(2n\log(k/m)) = (2/k)(1 - F_l(2n\log(k/m))). \tag{5}$$

We now consider some specific cases.

**Exponential** Suppose that all nodes have an equal and constant hazard rate, $1/n$ of being removed. Again, our analysis will not differentiate between the discrete time case where we remove a random node each time step; the discrete time case where each node's life span is a geometric random variable with mean $n$; and the continuous time case where each node's life span is an exponential random variable with mean $n$. Thus $F_l(t) = 1 - e^{-t}$ and by (5), $\Pr[k(t) = k] = 2m^2/k^3$ for $k \geq m$. The degree distribution is a power-law with the same exponent as a (static) preferential attachment network.

**Other exponents** The above analysis can be tweaked to obtain power-law degree distributions with other exponents. Specifically, $\Pr[k(t) = k] = \Theta(k^{-\gamma})$ for $\gamma > 1$ when $1 - F_l(t) = \exp(-t(\gamma - 1)/(2n) + o(t))$, that is, each node's hazard rate of removal is $(\gamma - 1)/(2n)$. When $\gamma = 3$ we can simply use an exponential distribution as discussed previously, but for other values of $\gamma$ we need to modify the exponential distribution to ensure that the total rate of node removal is 1 per unit time. When we have thinner tails, $\gamma > 3$, we can increase the hazard rates for small $t$ (i.e., for younger nodes) to compensate. For fatter tails, $1 < \gamma < 3$, we can compensate by either decreasing the hazard rates for small $t$ or truncating the degree distribution (i.e., setting a maximum node age).

**FIFO** Another special case is when we remove nodes in a first in first out basis, that is we remove the oldest node first. We can model this as a discrete time process where every step we first remove the oldest node and then add a new node. This is essentially the above case with $\gamma = 1$ and a truncated degree distribution. Thus the age of a random node is distributed $U[0, n]$. Hence by (5), $\Pr[k(t) = k] = 2/k$, for $k$ between $m$ and $m\sqrt{e}$.

## 5 Appendix

*Proof of Theorem 1.* From lemma 6, $TV(E_k(t), \bar{E}_k) = TV(X_k(p(t)), X_k(p))$ where $X_k(q) := \text{Binomial}(k, q)$ and $p(t) := p + (1-p)e^{-(\lambda+\mu)t}$. Define $A_k(q) := \{i : \Pr[X_k(p) = i] \geq \Pr[X_k(q) = i]\}$. When $q \geq p$, $A_k(q) = \mathbb{Z} \cap [0, ka_k(q)]$. From the pmf of the binomial distribution we can solve for $a_k$ and obtain $a_k(q) = \log \frac{1-q}{1-p} / \log \frac{p(1-q)}{q(1-p)}$. By lemma 7,

$$TV(X_k(p(t)), X_k(p)) = \Pr[X_k(p) \leq ka_k(p(t))] - \Pr[X_k(p(t)) \leq ka_k(p(t))]. \tag{6}$$

Focusing on the first claim, suppose that $t = \frac{1}{2(\lambda+\mu)}(\log k + b(p))$ as $k \to \infty$. Then $k^{1/2}(p(1-$

$p))^{-1/2}(a_k(p(t))-p) \to (b(p)/2)\sqrt{(1-p)/p}$ and similarly, $k^{1/2}(p(t)(1-p(t)))^{-1/2}(a_k(p(t))-p(t)) \to -(b(p)/2)\sqrt{(1-p)/p}$. Applying the CLT to (6),

$$TV(X_k(p(t)), X_k(p)) \to 2\left(\phi\left(\frac{b(p)}{2}\sqrt{(1-p)/p}\right) - \frac{1}{2}\right). \tag{7}$$

Thus by choosing $b(p)$ appropriately we prove the first claim..

Now focusing on the second claim, suppose that $t = \frac{\log(k/b(p))}{\lambda+\mu}$ as $k \to \infty$. Note that $\lambda + \mu = \frac{\alpha}{p(1-p)}$. Then $kp(t) \to c+b(p)$ and $ka_k(p(t)) \to b(p)/\log(1+b(p)/c)$. Now applying the law of *small numbers* to (6),

$$TV(X_k(p(t)), X_k(p)) \to TV(\text{Poisson}(c), \text{Poisson}(c+b(p))). \tag{8}$$

Thus by choosing $b(p)$ appropriately we prove the second claim. $\square$

**Lemma 6.** *Let $X_i \sim \text{Bernoulli}(p)$ and $Y_i \sim \text{Bernoulli}(q)$ be iid for $i = 1, \ldots, k$, respectively. Let $X^k := (X_1, \ldots, X_k)$ and $Y^k := (Y_1, \ldots, Y_k)$. Then $TV(X^k, Y^k) = TV(\text{Binomial}(k, p), \text{Binomial}(k, q))$.*

*Proof.* The sample space is the set of $k$-bit strings, $\Omega := \mathbb{Z}_2^k$. For $\omega := (\omega_1, \ldots, \omega_k) \in \Omega$, let $|\omega|$ denote the number of 1s in $\omega$. Using lemma 7,

$$TV(X^k, Y^k) = \frac{1}{2}\sum_\omega \left|p^{|\omega|}(1-p)^{k-|\omega|} - q^{|\omega|}(1-q)^{k-|\omega|}\right| \tag{9}$$

$$= \frac{1}{2}\sum_{i=0}^k \binom{k}{i}\left|p^i(1-p)^{k-i} - q^k(1-q)^{k-i}\right| \tag{10}$$

$$= TV(\text{Binomial}(k, p), \text{Binomial}(k, q)). \tag{11}$$

$\square$

**Lemma 7.** *Let $E := \{\omega : \xi(\omega) \geq \nu(\omega)\}$ and let $\bar{E}$ be its complement. Then $TV(\xi, \nu) = \xi(E) - \nu(E) = \nu(\bar{E}) - \xi(\bar{E})$. Further if the sample space is countable, then $TV(\xi, \nu) = (1/2)\sum_\omega |\xi(\omega) - \nu(\omega)|$.*

*Proof.* Note that $1 = \xi(E) - \nu(E) + \nu(E) + \xi(\bar{E}) = \nu(\bar{E}) - \xi(\bar{E}) + \nu(E) + \xi(\bar{E})$. Thus $\xi(E) - \nu(E) = \nu(\bar{E}) - \xi(\bar{E}) \geq 0$. From the construction of $E$, $TV(\xi, \nu) = \max\{\xi(E) - \nu(E), \nu(\bar{E}) - \xi(\bar{E})\}$, proving the first claim. Now $TV(\xi, \nu) = (1/2)((\xi(E) - \nu(E)) + (\nu(\bar{E}) - \xi(\bar{E})))$ proving the second claim. $\square$

*Proof of theorem 3.* Let $p = f(t) := 1 - e^{-\lambda t}$, $p_n := \log n / n$, and $t_n := \log n / (n\lambda)$. Choose an $\epsilon > 0$. Since edges are only added, the graph is connected at time $t$ iff $\tau_n \leq t$. At time $t$ the graph is a $G(n, p)$ graph with $p = f(t)$. Let $E_1$ be the event that a $G(n, p)$ graph with $p = p_1 := f((1-\epsilon)t_n) = (1 - \epsilon + o(1))p_n$ is connected. Define $E_2$ similarly with $p = p_2 := f((1+\epsilon)t_n) = (1 + \epsilon + o(1))p_n$. Note that $p_1 = (1 - \epsilon + o(1))p_n \leq (1 - \epsilon/2)p_n$ for large $n$ and similarly $p_2 \geq (1 + \epsilon/2)p_n$ for large $n$. Then by theorem 2, $\Pr[E_1] \to 0$ and $\Pr[E_2] \to 1$. Thus $\Pr[|\tau_n/t_n - 1| \leq \epsilon] \to 1$ proving the claim. □

*Proof of lemma 4.* The hitting time solves the recurrence,

$$((N - k)\lambda + k\mu + k\beta)t_k = 1 + (N - k)\lambda t_{k+1} + k\mu t_{k-1}. \tag{12}$$

Now define $x := kN^{-1/2}$ and $T(x) := t_k N^{1/2}$. Then this recurrence becomes

$$(-\lambda + x(\lambda + \mu)N^{-1/2})\frac{T(x + N^{-1/2}) - T(x)}{N^{-1/2}}$$
$$- (x\mu/N)\frac{T(x + N^{-1/2}) - 2T(x) + T(x - N^{-1/2})}{N^{-1}} + \beta x T(x) = 1. \tag{13}$$

As $N \to \infty$, this becomes the ODE, $-\lambda T'(x) + \beta x T(x) = 1$ whose solution is

$$T(x) = e^{x^2\beta/(2\lambda)}\left(T(0) - \frac{1}{\lambda}\int_0^x e^{-s^2\beta/(2\lambda)}ds\right). \tag{14}$$

Since $t_k$ is decreasing and nonnegative, $T(x)$ must also be. Thus

$$T(0) = \frac{1}{\lambda}\int_0^\infty e^{-s^2\beta/(2\lambda)}ds = \sqrt{\pi/(2\beta\lambda)}, \tag{15}$$

proving the claim.

□

*Proof of theorem 5.* We determine a lower bound by constructing a graph $G(t)$ where edges are never removed and looking at the time until it is connected. Unlike the case of $\beta = \infty$, we assume that the time until an edge appears between a pair of nodes is $Z := \text{Exp}(\lambda) + \text{Exp}(\beta)$, the time it takes for an edge to appear in the original graph and then the infection to travel across it.

Then $G(t)$ equals in distribution a $G(n,p)$ graph with $p = f(t) := \Pr[Z \leq t]$. We then take the same approach as in theorem 3 with the modified $f(t)$ and $t_n := \sqrt{2\log n/(\beta\lambda n)}$. Note that $f(t) = 1 - (\lambda e^{-\beta t} - \beta e^{-\lambda t})/(\lambda - \beta)$, $f((1+\epsilon)t_n) = ((1+\epsilon)^2 + o(1))p_n \geq (1+\epsilon)p_n$ for large $n$, and $f((1-\epsilon)t_n) = ((1-\epsilon)^2 + o(1))p_n \leq (1-\epsilon)p_n$ for large $n$. $\square$

# References

[1] A.-L. Barabási and R. Albert. Emergence of Scaling in Random Networks. *Science*, 286(5439):509–512, 1999. URL http://www.sciencemag.org/cgi/doi/10.1126/science.286.5439.509.

[2] A.-L. Barabási and R. Albert. Emergence of scaling in random networks. *Science*, 286(5439):509–512, 1999.

[3] A. Barrat, M. Barthélemy, and A. Vespignani. *Dynamical Processes on Complex Networks*. Cambridge University Press, 2008.

[4] A. Bonato. A survey of models of the web graph. In A. Lpez-Ortiz and A. Hamel, editors, *Combinatorial and Algorithmic Aspects of Networking*, volume 3405 of *Lecture Notes in Computer Science*, pages 159–172. Springer Berlin / Heidelberg, 2005.

[5] A. Cami and N. Deo. Techniques for analyzing dynamic random graph models of web-like networks: An overview. *Netw.*, 51:211–255, 2008. URL http://portal.acm.org/citation.cfm?id=1378692.1378693.

[6] F. R. K. Chung and L. Lu. Coupling online and offline analyses for random power law graphs. *Internet Mathematics*, 1(4), 2003.

[7] C. Cooper, A. M. Frieze, and J. Vera. Random deletion in a scale-free random graph process. *Internet Mathematics*, 1(4), 2003.

[8] K. Dietz and K. P. Hadeler. Epidemiological models for sexually transmitted diseases. *Journal of Mathematical Biology*, 26(1):1–25, 1988. URL http://www.ncbi.nlm.nih.gov/pubmed/3351391.

[9] Durrett. *Random Graph Dynamics*. Cambridge University Press, 2007.

[10] P. Erdős and A. Rényi. On random graphs. *Publicationes Mathematicae*, 6(26):290–297, 1959.

[11] N. H. Fefferman and K. L. Ng. How disease models in static networks can fail to approximate disease in dynamic networks. *Phys. Rev. E*, 76(3):031919, 2007.

[12] J. P. Ferry, D. Lo, S. T. Ahearn, and A. M. Phillips. Network detection theory. In N. Memon, J. D. Farley, D. L. Hicks, and T. Rosenorn, editors, *Mathematical Methods in Counterterrorism*, pages 161–185. Springer, 2009.

[13] E. N. Gilbert. Random graphs. *Ann. Math. Statist.*, 30:1141–1144, 1959.

[14] T. Gross and B. Blasius. Adaptive coevolutionary networks: a review. *J. R. Soc. Interface*, 5(20):259–271, 2008.

[15] T. Gross, C. J. D. D'Lima, and B. Blasius. Epidemic dynamics on an adaptive network. *Phys. Rev. Lett.*, 96(20):208701, 2006.

[16] P. Holland and S. Leinhardt. A dynamic model for social networks. *The Journal of Mathematical Sociology*, 5(1):5–20, 1977.

[17] M. Huisman and T. Snijders. Statistical analysis of longitudinal network data with changing composition. *Sociological Methods & Research*, 32(2):253, 2003.

[18] C. Kamp. Demographic and behavioural change during epidemics. *Procedia Computer Science*, 1(1):2253 – 2259, 2010. URL http://www.sciencedirect.com/science/article/B9865-506HM1Y-8Y/2/526ba0b24caf5da48f456da884474ba4. ICCS 2010.

[19] C. Kamp. Untangling the interplay between epidemic spread and transmission network dynamics. *PLoS Comput Biol*, 6(11):e1000984, 2010. URL http://dx.doi.org/10.1371/journal.pcbi.1000984.

[20] J. K. Ochab and P. F. Goára. Shift of percolation thresholds for epidemic spread between static and dynamic small-world networks. 2010. URL http://arxiv.org/abs/1011.2985v1.

[21] R. Pastor-Satorras and A. Vespignani. Epidemic spreading in scale-free networks. *Physical Review Letters*, 86:3200–3203, 2001.

[22] Y. Peres, A. Sinclair, P. Sousi, and A. Stauffer. Mobile geometric graphs: Detection, coverage and percolation. In *Proceedings of the ACM-SIAM Symposium on Discrete Algorithms 2011*. 2011.

[23] S. Risau-Gusman. Influence of network dynamics on the spread of sexually transmitted diseases. 2010. URL `http://arxiv.org/abs/1004.1378v1`.

[24] I. B. Schwartz and L. B. Shaw. Rewiring for adaptation. *Physics*, 3:17, 2010.

[25] Y. Schwarzkopf, A. Rákos, and D. Mukamel. Epidemic spreading in evolving networks. *Phys. Rev. E*, 82(3):036112, 2010.

[26] L. B. Shaw and I. B. Schwartz. Fluctuating epidemics on adaptive networks. *Phys. Rev. E*, 77(6):066101, 2008.

[27] T. Snijders. Models for longitudinal network data. In *Models and methods in social network analysis*, chapter 11, pages 215–247. Cambridge Univ Press, 2005.

[28] T. A. Snijders and P. Doreian. Introduction to the special issue on network dynamics. *Social Networks*, 32(1):1–3, 2010. URL `http://www.sciencedirect.com/science/article/B6VD1-4Y5BDH7-2/2/566bc1f45c261c58f6194749ae37aa96`. Dynamics of Social Networks.

[29] S. Van Segbroeck, F. C. Santos, and J. M. Pacheco. Adaptive contact networks change effective disease infectiousness and dynamics. *PLoS Comput Biol*, 6(8):e1000895, 2010. URL `http://dx.doi.org/10.1371%2Fjournal.pcbi.1000895`.

[30] E. Volz and L. A. Meyers. Susceptible infected recovered epidemics in dynamic contact networks. *Proceedings of the Royal Society B: Biological Sciences*, 274(1628):2925–2934, 2007. URL `http://rspb.royalsocietypublishing.org/content/274/1628/2925.abstract`.

[31] S. Wasserman. *Stochastic Models for Directed Graphs*. Ph.D. thesis, Harvard, 1977.